\documentclass[12pt]{article}

\usepackage{epsfig}

\def\absname{ABSTRACT:} 
\def\abs{\if@twocolumn
\section*{\absname}
\else \small
\noindent {\bf \absname}
\fi}
\def\endabs{\if@twocolumn\else\endquotation\fi}
\mark{{}{}}

\def\resumoname{RESUMO:} 
\def\resumo{\if@twocolumn
\section*{\resumoname}
\else \small
\noindent {\bf \resumoname}
\fi}
\def\endresumo{\if@twocolumn\else\endquotation\fi}
\mark{{}{}}

\font\tenmsbm=msbm10 scaled \magstep 1
\font\sevenmsbm=msbm10
\font\fivemsbm=msbm5
\newfam\msbmfam
\textfont\msbmfam=\tenmsbm
\scriptfont\msbmfam=\sevenmsbm
\scriptscriptfont\msbmfam=\fivemsbm
\def\msbm{\fam\msbmfam}

\newcommand{\C}{{\msbm C}}

\def\halmos{\hfill\rule{6pt}{6pt}}

\newtheorem{lemma}{\indent Proposition}[section]
\newtheorem{lem}[lemma]{\indent Lemma}
\newtheorem{corol}[lemma]{\indent Corollary}

\title{Spectral Transformation Algorithms for Computing Unstable Modes
of Large Scale Power Systems} 

\author{L. H. Bezerra  \\
{\small Departamento de Matem\'atica} \\
{\small Universidade Federal de Santa Catarina} \\
{\small Florian\'opolis, 88040-900, Brazil} \\
{\small e-mail: licio@mtm.ufsc.br} 
\and C. Tomei \\
{\small Departamento de Matem\'atica} \\
{\small Pontif\'{\i}cia Universidade Cat\'olica} \\
{\small Rio de Janeiro, 22453-900, Brazil} \\
{\small e-mail: carlos@mat.puc-rio.br}}

\date{}

\begin{document}

\maketitle


\begin{abs}
In this paper we describe spectral transformation algorithms
for the computation of eigenvalues with positive
real part of sparse nonsymmetric matrix pencils $(J,L)$, where $L$
is of the form $\pmatrix{M&0\cr 0&0}$.
For this we define a different extension of M\"obius transforms
to pencils that inhibits the effect on iterations of the spurious eigenvalue
at infinity.
These algorithms use a technique of preconditioning the initial vectors
by M\"obius transforms which together with shift-invert iterations
accelerate the convergence to the desired eigenvalues.
Also, we see that M\"obius transforms can be successfully used in inhibiting
the convergence to a known eigenvalue. Moreover, the procedure
has a computational cost similar to power or shift-invert iterations with
M\"obius transforms: neither is more expensive than the usual
shift-invert iterations with pencils. Results from tests with
a concrete transient
stability model of an interconnected power system whose
Jacobian matrix has order 3156 are also reported here.
\end{abs}


\vskip 1 pc

\begin{resumo}
Neste artigo, descrevemos algoritmos baseados em transforma\c c\~oes
espectrais para computa\c c\~ao de autovalores com parte real positiva
de {\it pencils} de matrizes esparsas e n\~ao sim\'etricas,
$(J,L)$, em que $L$ \'e da forma $\pmatrix{M&0\cr 0&0}$.
Para isso definimos uma extens\~ao das transformac\~oes de M\"obius a
{\it pencils} que inibe a atua\c c\~ao do autovalor infinito
sobre as itera\c c\~oes.
Esses algoritmos usam uma t\'ecnica de precondicionamento dos vetores
iniciais via transformadas de M\"obius que junto com itera\c c\~oes tipo
pot\^encia inversa com {\it shift}
aceleram a converg\^encia para os autovalores desejados.
Vemos tamb\'em que as transformadas de M\"obius podem ser usadas com sucesso
no processo de inibir a converg\^encia para um autovalor j\'a conhecido.
Al\'em disso, esse procedimento tem um custo computacional semelhante ao
custo computacional de itera\c c\~oes tipo pot\^encia ou pot\^encia inversa
com {\it shift}: t\~ao caro como itera\c c\~oes tipo pot\^encia inversa com
{\it shift} aplicadas em {\it pencils}.
S\~ao tamb\'em apresentados aqui resultados de testes com
um modelo pr\'atico para o problema de
estabilidade transiente de um sistema de pot\^encia interconectado,
cuja matriz jacobiana \'e de ordem 3156.
\end{resumo}

\noindent {\small {\bf PALAVRAS-CHAVE:} autovalores, estabilidade,
transforma\c c\~oes de M\"obius,
autovalores generalizados}

\vskip 1 pc



\section{Introduction}

The power system eletromechanical stability problem can be described
by a nonlinear system of differential and algebraic
equations
$$
\left\{ \begin{array}{ll}
               \dot x = f(x,y) \\
                    0 = g(x,y),    \end{array}
\right. \eqno(1.1)
$$
where x, the state vector, contains the dynamic variables and y,
the algebraic variables.
After linearization around a system operating point ($x_0,y_0$), i.e,
($x_0,y_0$) such that f($x_0,y_0$) = 0, equation (1.1) becomes
$$
\pmatrix{\Delta \dot x \cr 0} = \pmatrix{J_1&J_2 \cr J_3&J_4}
\pmatrix{\Delta x \cr \Delta y}. \eqno(1.2)
$$
By eliminating the vector $\Delta y$ in (1.2) we obtain
$$
\Delta \dot x = A \Delta x,  \eqno (1.3)
$$
where $ A = J_1 - J_2 J_4^{-1} J_3 $ represents the system state matrix,
whose eigenvalues provide information about the singular point
local stability
of the non-linear system. The symbol $\Delta$ used to represent
an incremental change from a steady-state value will be omitted
from now on.

By a classical result of ODE theory, the local stability of the system
(1.1) can be predicted from the system (1.3).
If $A$ is diagonalizable the solution of (1.3) is a sum of vectors of the
type $e^{\lambda_i t}v_i$, where $v_i$ is an eigenvector associated with the
eigenvalue $\lambda_i$. Thus, if all eigenvalues have negative real part,
the solution decays to zero and in this case, the system is called stable.
For an eigenvalue with positive real part, the absolute value of this expression
increases in time and the system is unstable ---
these eigenvalues are called unstable modes.
Eigenvalues with null real part give rise to
oscillation, which never disappears. Also,
eigenvalues with negative real part and
non-zero imaginary part, but with
small ratio between the real and imaginary parts,
cause an oscillation which takes a long time to disappear ---
these are the low damped modes of the system.
In the power system stability problem, we consider a mode $\lambda$ to
be low damped if $|Re \, \lambda|<0.02 |Im \, \lambda|$.
In this paper, we search for algorithms which solve the local stability
problem (1.1) by computing eigenvalues:
we search for
the eigenvalues of $A$ with positive real part.
The state matrices $A$ are real, non-symmetric and  dense, usually
too large for the computation of eigenvalues by the QR method.
On the other hand, the Jacobian matrices $J$ are sparse and
linear systems with $J$  may be solved by variants of Gaussian elimination.
The eigenvalue problem for A,
$$
A x = \lambda x,  \eqno(1.4)
$$
can be stated equivalently in terms of the Jacobian matrix
$J=\pmatrix{J_1&J_2 \cr J_3&J_4}$
so that
$$
\pmatrix{J_1&J_2 \cr J_3&J_4} \pmatrix{x \cr y} =
\lambda \pmatrix{x \cr 0}.  \eqno(1.5)
$$
or, in matrix notation,
$$
J z = \lambda L z, \eqno(1.6)
$$
where $L=\pmatrix{I&0\cr 0&0}$ is
the singular diagonal matrix with 1 (resp. 0)
at diagonal entries related to state (resp. algebraic) variables and
$z=\pmatrix{x \cr y}$ is an eigenvector of
$(J-\lambda L)$. Since $z\ne 0$, $\lambda$ is said to be an eigenvalue of
the pencil $(J,L)=\{ J-aL|a\in \C\}$.
A possible approach to search for unstable modes
is to use the shift-invert transform
$$ (J- \lambda_k L) z_{k+1} = L z_k, \eqno(1.7) $$
with initial shifts on the imaginary axis \cite{nm}.
Although in \cite{nm}
the problem is not described as the generalized
eigenvalue problem (1.6), the system (1.7) is solved in order to
implicitly calculate $(A - \lambda I)^{-1} x_k$ --- in the authors
words, they make use of the augmented system associated to $J$
to shift-invert the state matrix $A$.
Other methods like subspace iteration
and Arnoldi methods were also adapted to augmented systems
and some results of their applicability in power systems are reported
in \cite{ws}.
The use of the Cayley transform technique
in order to find rightmost eigenvalues
of a non-symmetric matrix was probably first reported
in 1987 at the IEEE PICA Conference. There, Uchida and Nagao
proposed to search for the biggest eigenvalues in absolute value
of $S=(A+hI)(A-hI)^{-1}$, $h>0$
to detect unstable modes of the state matrix
$A$ of a power system \cite{jp}.
There, however, the matrix-vector multiplication
$Sv$ was performed in a rather innefficient way.
In \cite{lh}, this operation was better implemented and the
Cayley transform was extended in two different ways, described in
detail in \S 2, in order to solve
the power system stability problem, given as
a generalized eigenvalue problem $Jz=\lambda Lz$,
where $J$ is not symmetric and $L$ is diagonal with elements 1 or 0
along the diagonal.
The use of the Cayley transform technique
to find rightmost eigenvalues of the problem $Jx=\lambda Lx$,
for non-symmetric $J$
and $L=\pmatrix{M&0\cr 0&0}$, is also
found in Computational Fluid Dynamics \cite{cg}, \cite{ga}; an overview
of this technique in several areas is \cite{me};
also, in ARPACK \cite{le}, Arnoldi iterations can be performed with
Cayley transforms.
All these references make use of the extension
$(J,L) \mapsto (J-\sigma L)^{-1}(J+\overline{\sigma} L)$, which is one
of the two extensions of the Cayley transform analysed in \cite{lh},
which turns out not to be the best for the problem of interest.
Indeed, this extension requires a non-obvious
strategy to control instability caused by the spurious eigenvalue
at infinity of the generalized eigenvalue problem,
in the case when $L$ is singular.
Here we propose to consider another extension:
$(J,L) \mapsto (J+\overline{\sigma}L)(J-\sigma L)^{-1}$.
For this matrix the eigenspace associated with the eigenvalues that
correspond to finite eigenvalues of the original problem is just
the range of $L$, as we shall see in \S 2.
Thus, the spurious eigenvalue is handled by keeping iterations
in this space.
In the power system stability problem, the range of $L$ is the
set of vectors with coordinates related to algebraic variables equal to zero.
In the case $L=\pmatrix{M&0\cr 0&0}$ is a matrix of order $n+m$, where
$M$ is a $n\times n$ positive definite matrix, as in \cite{cg},
the space is just the set of vectors with
the last m coordinates equal to zero.
In \S 2, we introduce M\"obius transforms (a generalization of the
Cayley transforms) for the generalized non-symmetric eigenvalue problem.
M\"obius transforms can be used to precondition random initial vectors
as well as to
inhibit the convergence to eigenvalues already found
(much like a deflation technique), at a computational cost
not more expensive than a shift-invert iteration with pencils.
The application of polynomial filters to vectors in the computation of
eigenvalues of sparse non-symmetric matrices has been the subject of several
papers \cite{so}, \cite{sa}, \cite{ch}.
We will see that M\"obius transforms can be seen as the action of a
rational function filter
which gives infinite and zero weights to two arbitrary
points in the complex plane.
These techniques are presented in \S 3
together with two algorithms based on them.
Finally, in \S 4
we apply four methods to the computation of the unstable modes of
a pencil $(J,L)$, where $J$ is a sparse matrix of order 3156
and $L$ is diagonal with elements 1 or 0 on the diagonal, of rank 790:
the algorithms we suggest, an implementation of the Arnoldi method (ARPACK) and
an implementation of the Arnoldi method with acceleration by Chebyshev
polynomials (ARNCHEB), both applied on the extension of a Cayley transform
introduced here
and a subspace iteration applied on the pencil $(J,L)$.


\section{The Generalized Eigenvalue Problem and M\"obius Transforms}

The M\"obius transforms are complex functions
\begin{eqnarray*}
c_{k,\alpha,\beta}:\C\cup\{\infty\} &\longrightarrow &
\C\cup\{\infty\} \\
\mbox{$s$\hspace{22pt}}&\mapsto & k\,\frac{s+\overline{\beta}}{s-\alpha}
\end{eqnarray*}
where $\alpha + \overline{\beta} \ne 0$.
They are conformal mappings which map lines and circumferences to
lines or circumferences.
When $k=1$ and $\beta=\alpha$, with $Re \, \alpha \ne 0$, these functions
are the so called Cayley transforms,
which map the semiplane $\{ z | Re\, z>0\}$ to $\{ z ; |z|>1\}$
($\{ z; |z|<1\}$) if $Re \, \alpha >0$ ($Re \, \alpha <0$).
M\"obius transforms can be defined in the space of square matrices in a
analogous way: given a square matrix $A$ and $\alpha$ not in
$\lambda (A)$ (the spectrum of $A$) then
$k (A+\overline{\beta}I)(A-\alpha I)^{-1}$ is a matrix whose spectrum is
$\{ k(\lambda + \overline{\beta})(\lambda - \alpha)^{-1};
\lambda \in \lambda(A) \}$.
However, the extension of these transforms to pencils $(J,L)$
can be done in two ways.
The goal of this section is to present the advantages of one extension
over the other in the application to
the power system stability problem. Also, we will see that
multiplication of a M\"obius transform against a vector
requires no more computations as solving
the equation $(J-aL)z=w$, for some (easily computed) scalar $a$.
We begin with a lemma that insures that under generic conditions the
pencil $(J,L)$ has exactly $n$ eigenvalues, where $n$ is the rank of $L$.
Since by the Singular Value Decomposition there are unitary matrices $U,V$
such that $U(J-aL)V^H=\hat{J}-a\Sigma$,
$\Sigma=diag(d_1,...,d_n,0,...,0)$ for appropriate real numbers $d_i$
\cite{go}, we may suppose without loss
that $L=diag(d_1,...,d_n,0,...,0)$.

\begin{lem}
Let $J=\pmatrix{J_1&J_2\cr J_3&J_4}$ and
$L=\pmatrix{D&0\cr 0&0}$, where $D$ is a nonsingular diagonal matrix. Then,
if $J_4$ is nonsingular,
$(J-\sigma L)$ is singular if and only if $\sigma$ is an eigenvalue of
$D^{-1}A$, where $A=J_1-J_2J_4^{-1}J_3$.
\end{lem}

\indent {\em Proof:}
If $J_4$ is nonsingular,
$A=J_1 - J_2 J_4^{-1} J_3$ can be defined. Now,
$$ \pmatrix{J_1 - \sigma D & J_2 \cr J_3 & J_4} \pmatrix{x\cr
y}=\pmatrix{0\cr 0} \Longleftrightarrow
(A - \sigma D)\, x = 0 \mbox{ and } y=-J_4^{-1} J_3 \, x.$$ Therefore,
$(J-\sigma L)$ is singular $\Longleftrightarrow$ $\sigma$ is
an eigenvalue of $D^{-1}A$.
\halmos

\vskip 1pc

From now on, the non-sigularity of $J_4$ will be assumed.
Let $\alpha \in \C$, such that $(J-\alpha L)$ is nonsingular.
For $k,\beta \in \C$, $k\ne 0$ and $\beta \ne -\bar{\alpha}$, let
$$ C_{k,\alpha,\beta} = k(J+\bar{\beta} L)(J-\alpha L)^{-1}
\mbox{ and } D_{k,\alpha,\beta} = k(J-\alpha L)^{-1}(J+\bar{\beta} L).
$$
Notice that these two matrices have the same spectrum.
Moreover, they are both extensions of M\"obius transforms
applied to the pencil $(J,L)$ and the spectrum
of $(J,L)$ is related to the spectrum of the two extensions
by the relation
$s\mapsto k\displaystyle{\frac{s+\overline{\beta}}{s-\alpha}}.$
However, the eigenspaces of both extensions are not the same in general,
according to the
following propositions whose proofs are left to the reader.

\begin{lemma}
(a) $z\ne 0$ is an eigenvector of $C_{k,\alpha,\beta}$ associated with $\mu$,
$\mu \ne k$, if and only if
$Jw = \lambda L w$, where
$w=(J-\alpha L)^{-1} \, z$ and
$\lambda =\displaystyle{\frac{k\bar{\beta} + \mu \, \alpha}{\mu -k}}$.

(b) $z\ne 0$ is an eigenvector of $C_{k,\alpha,\beta}$ associated with
the eigenvalue $k$ if and only if $L (J-\alpha L)^{-1}z=0$.
\end{lemma}

\begin{lemma}
(a) $z\ne 0$ is an eigenvector of $D_{k,\alpha,\beta}$ associated with $\mu$,
$\mu \ne k$, if and only if
$Jz = \lambda L z$, where
$\lambda =\displaystyle{\frac{k\bar{\beta} + \mu \, \alpha}{\mu -k}}$.

(b) $z\ne 0$ is an eigenvector of $D_{k,\alpha,\beta}$ associated with
the eigenvalue $k$ if and only if $Lz=0$.
\end{lemma}

Thus, the finite eigenvalues of the pencil $(J,L)$ correspond
to the eigenvalues of $C_{k,\alpha,\beta}$ or $D_{k,\alpha,\beta}$
that are different from $k$.
These eigenspaces are described in the following proposition.

\begin{lemma}
(a) The eigenspace of $C_{k,\alpha,\beta}$ associated with the eigenvalues
different from $k$ is the range of $L$.

(b) The eigenspace of $D_{k,\alpha,\beta}$ associated with the eigenvalues
different from $k$ is the range of $(J-\alpha L)^{-1} L$.
\end{lemma}
\vskip -6 pt

{\em Proof:}
Let $z$ be a vector such that $C_{k,\alpha,\beta}^H \, z = k z$. Thus
$(J^H-\beta L^H) z =(J^H-\bar{\alpha}L^H) z$ and,
since $\beta\ne \bar{\alpha}$, $L^H z=0$.
Since the eigenspace of $C_{k,\alpha,\beta}$ associated with the eigenvalues
different from $k$ is orthogonal to
the eigenspace of $C_{k,\alpha,\beta}^H$ associated with $k$,
that eigenspace is the range of $L$.

The proof of the second part of the proposition is analogous.
\halmos

\begin{corol}
Let $J=\pmatrix{J_1&J_2\cr J_3&J_4}$ a matrix of order $n+m$, where $J_4$ is
a $m\times m$ nonsingular matrix, and let
$L=\pmatrix{M&0\cr 0&0}$, where $M$ is a $n\times n$ nonsingular matrix.
Then the eigenspace of $C_{k,\alpha,\beta}$ associated with the eigenvalues
different from $k$ is the space of vectors
whose last $m$ coordinates are zero.
\end{corol}

Since we are interested in calculating eigenvalues of the pencil $(J,L)$
from M\"obius transforms, the eigenvalue $k$ of the transform must be treated
with special care.
The corollary above states that in the case of the power system
stability problem, for instance,
the invariant subspace $V$ corresponding to the eigenvalues different from
$k$ is the set of vectors $v$ that have null coordinates
in the positions related to algebraic variables. The iterates
of our approximate eigenvectors ought to stay in this subspace, and
if the computation of $C_{k,\alpha,\beta}v$ leaves $V$
because of errors due to finite precision arithmetic,
we simply project the results back to $V$ by zeroing the appropriate
coordinates. The analogous iteration
with the M\"obius extension $D_{k,\alpha,\beta}$ is not
subject to such an easy stabilization procedure, and in
this case approximate eigenvectors will frequently
converge to the eigenspace associated to $k$, which is of
no real interest for the pencil eigenvalue problem.
Hence, we will only consider here iterations making use of
the extension $C_{k,\alpha,\beta}$.

Now, let $\sigma \in \C$, with $Re \,\sigma > 0$, such that $(J-\sigma L)$
is nonsingular, and consider
$$ C_{\sigma} = C_{1,\sigma,\sigma} =
(J+\overline{\sigma} L)(J-\sigma L)^{-1} \eqno(2.1) $$
Let $\mu\ne 1$. Then, simple calculations obtain
$$
(C_{\sigma} - \mu I) =
(1-\mu)I + 2 \, Re \,\sigma \,L \, (J - \sigma \, L)^{-1}
$$
$$
(C_{\sigma} - \mu I)^{-1} = \frac{1}{1-\mu}
\; [\, I + \frac{2 \, Re \,\sigma}{\mu -1} \,L \,
(J-\frac{\bar{\sigma} + \mu \sigma}{\mu -1}\, L)^{-1} \,],
$$
$$
(C_{\sigma}^{-1} - \mu I)^{-1} = \frac{1}{1-\mu}
\; [\, I - \frac{2 \, Re \,\sigma}{\mu -1} \,L \,
(J+\frac{\sigma + \mu \bar{\sigma}}{\mu -1}\, L)^{-1} \,].
$$
Thus, the computational cost of applying
any of the three matrices in the left-hand side to a vector
is equivalent to solving a system $(J-a L)w=z$.
Similar results hold for both extensions of M\"obius transforms.

The spectra of the pencil $(J,L)$ and of $C_{\sigma}$
are related to each other by the bijective function
\begin{eqnarray*}
c_{\sigma}:\C\cup\{\infty\} &\longrightarrow &
\C\cup\{\infty\} \\
\mbox{$s$\hspace{22pt}}&\mapsto & \frac{s+\overline{\sigma}}{s-\sigma}
\end{eqnarray*}
This function maps complex numbers with negative real part onto the
unitary circle, pure imaginary ones onto the unitary
circumference and those with positive real part onto numbers
outside the unitary circle. The search for eigenvalues of largest
absolute value of the extension of the Cayley transform thus
obtains the unstable modes of the original pencil. Unfortunately,
this transformation clusters
some eigenvalues very close together, thus affecting adversely
the rate of convergence.
Therefore, we need techniques to accelerate the convergence to
the desired eigenvalues. This is the goal of the next section.

\section{A Class of Spectral Algorithms}

Two ways of extending M\"obius transforms to the pencil $(J,L)$
were described in the previous section.
When $L$ is a matrix of order $n+m$ and
of the type $\pmatrix{M&0\cr 0&0}$, where
$M$ is an $n\times n$ nonsingular matrix, the eigenvectors of
the extension $C_{k,\alpha,\beta}$
associated with eigenvalues different from $k$ are the vectors which
have zero entries in the last $m$ coordinates.
Therefore, in this case,
the attraction of the eigenvalue $k$, which
corresponds to the infinite eigenvalue of the pencil, can be easily avoided.
Now, the power system stability problem,
where $M$ is the $n\times n$ identity matrix,
has additional features that lead us to explore M\"obius techniques.
Usually, this problem has several negative
real eigenvalues with large modulii, which are mapped to eigenvalues close
to $k$: the clustering of eigenvalues substantially reduces the speed
of convergence of the power method \cite{go}.
In this section we introduce a class of algorithms which use
M\"obius transforms to precondition vectors,
in order to yield more convenient initial vectors
for a search process with shift-invert
M\"obius transforms iterations inside the unit circle.
Also, we introduce a way of inhibiting known eigenvalues in the iteration
by yet another use of M\"obius transforms.

\subsection{Preconditioning and Shifts}

Let $\sigma \in \C$, $Re\, \sigma >0$.
As seen in \S 2, the eigenvalues of $(J,L)$ with positive real part
correspond to the eigenvalues of
$C_{\sigma}^{-1}$ located inside the unit circle.
In order to achieve larger components of eigenvectors of $C_{\sigma}^{-1}$
associated with eigenvalues of modulus less than
one in the search vectors,
we start with random vectors and apply
$C_{\sigma}$ to them a few times, obtaining the so called
preconditioned vectors.
We then use these vectors in a search process for eigenvalues
inside the unit circle. How should one choose $\sigma$ ?
The reality of the pencil $(J,L)$ implies that eigenvalues come in
conjugate pairs, which is still true for the matrix $C_{\sigma}$ if
$\sigma$ is taken to be real: our search for eigenvalues is then
reduced to, say, the upper half-disk.
Also, in this case, a simple differentiation shows that
$$\sigma = |\lambda|$$ maximizes
$|c_{\sigma}(\lambda)| = |\frac{\lambda + \sigma}{\lambda - \sigma}|$,
where $\lambda= a + b\, i$, $b\neq 0$ --- thus, the choice
$\sigma = |\lambda|$ takes eigenvalues of the pencil of absolute value
$|\lambda|$ to eigenvalues of $C_{\sigma}$ of largest possible
absolute value.
We then take shift-invert iterations with $C_{\sigma}^{-1}$ as the matrix
to be shifted.
Since the  eigenvalues of $C_{\sigma}^{-1}$ are the inverse of the ones
of $C_{\sigma}$,
the dominant eigenvectors of $C_{\sigma}$ are associated with
the eigenvalues of $C_{\sigma}^{-1}$ which are either inside the unit disk
or near the unit circle.
Based on these remarks we introduce the first algorithm.

\begin{itemize}
\item Algorithm I
\begin{itemize}
\item[Step 1]  Begin with $r$ orthogonal vectors belonging to
${\cal N}(L)^{\perp}$.
\item[Step 2] Multiply the vectors by $C_{\sigma}$ $p$ times, normalizing them
after each multiplication (e.g., keep them with sup norm equal to one).
Let $v_i$, $i=1,...,r$, be the resulting vectors.
\item[Step 3] Take initial shifts $\mu_k^{(0)}$, $k=1,...,s$, in a
circumference of radius $\epsilon$, $0< \epsilon \leq 1$.
\item[Step 4] Let $w_i^{(0)} = v_i$, $i=1,...,r$.

For $k=1,...,s$
    \begin{itemize}
    \item[] for $j=1,2,...$
            \begin{itemize}
            \item[] $u_i^{(j)} =
(C_{\sigma}^{-1} - \mu_k^{(j-1)} I)^{-1} w_i^{(j-1)}$;
            \item[] $w_i^{(j)} = u_i^{(j)} / \alpha_i^{(j)}$,
where $\alpha_i^{(j)}$ is the coordinate of $u_i^{(j)}$ of maximum absolute
value;
            \item[] if $j=t,2t,...$
                    \begin{enumerate}
                    \item[] $\mu_k^{(j)}=\mu_k^{(j-1)} +
                        \displaystyle{\frac{1}{\alpha_q^{(j)}}}$ ,
\ where $q$ is such that

$||w_q^{(j)} - w_q^{(j-1)}||_{\infty}$ =
$\displaystyle{\min_{1\leq i\leq r}}$
$||w_i^{(j)} - w_i^{(j-1)}||_{\infty}$.
                  \end{enumerate}
            \end{itemize}
    \end{itemize}
\end{itemize}

\end{itemize}

Preconditioning is performed in Step 2 ---
different choices of
$p$ are presented in the experiments of \S 4.
The shifts in Step 3 are taken in a circle centered at the origin.
The convergence of shift-invert iterations to an eigenvalue
depends enormously on the choice of shift:
the aim of this shifting
strategy is to cover the unit circle. It is at this point that the choice of
a real parameter $\sigma$ entitles us to divide by two the search for
eigenvalues by taking into account the reality of the original pencil.
Every $t$ iterations, the shift $\mu_k^{(j)}$ is updated by
the formula above \cite{wa}.
When the variation of one of the
vectors between
the previous and the current iterations is smaller than a
fixed tolerance, the shift is updated and the resulting value
is taken to be an eigenvalue.

\subsection{Inhibiting Convergence of Eigenvalues}

M\"obius transforms could be used to inhibit the convergence to an
eigenvalue $\xi$ of $C_{\sigma}^{-1}$ if in Algorithm I
iteration vectors were multiplied by $(C_{\sigma}^{-1} - \xi I)$.
Notice that after $t$ iterations with the matrix
$
(C_{\sigma}^{-1} - \xi I)(C_{\sigma}^{-1} - \mu_k^{(j)}I)^{-1}
$
(again a M\"obius transform)
the updated shift should be
$$\frac{\mu_k^{(j)} \alpha_q^{(j)} - \xi}{\alpha_q^{(j)} -1},$$
where $\alpha_q^{(j)}$ is defined in a similar way.
The computational cost is again no more expensive than
a standard shift-invert step for the generalized eigenproblem because
$$
(C_{\sigma}^{-1} - \xi I)(C_{\sigma}^{-1} - \mu_k^{(j)}I)^{-1}=
I + (\mu_k^{(j)} - \xi)(C_{\sigma}^{-1} - \mu_k^{(j)}I)^{-1}
$$
The algorithm below  uses this deflation-type strategy.

\begin{itemize}

\item Algorithm II

Use the Algorithm I to obtain  preconditioned vectors
$w_i^{(0)} = v_i$, $i=1,...,r$, and initial shifts
$\mu_k^{(0)}$, $k=1,...,s$.

For $k=1$
\begin{itemize}
\item[] Follow Algorithm I up to  convergence to an eigenvalue $\xi$.
\end{itemize}

For $k=2,...,s$
    \begin{itemize}
    \item[] for $j=1,2,...$
            \begin{itemize}
            \item[] $u_i^{(j)} =  (C_{\sigma}^{-1} -\xi I)
(C_{\sigma}^{-1} - \mu_k^{(j-1)} I)^{-1} w_i^{(j-1)}$;
            \item[] $w_i^{(j)} = u_i^{(j)} / \alpha_i^{(j)}$,
where $\alpha_i^{(j)}$ is the coordinate of $u_i^{(j)}$ of maximum absolute
value;
            \item[] if $j=t,2t,...$
                    \begin{enumerate}
                    \item[] $\mu_k^{(j)}=
 \displaystyle{\frac{\mu_k^{(j)} \alpha_q^{(j)} - \xi}{\alpha_q^{(j)} -1}}$ ,
\ where $q$ is such that

$||w_q^{(j)} - w_q^{(j-1)}||_{\infty}$ =
$\displaystyle{\min_{1\leq i\leq r}}$
$||w_i^{(j)} - w_i^{(j-1)}||_{\infty}$.
                  \end{enumerate}
            \end{itemize}
    \end{itemize}

\end{itemize}

\section{Tests and Comparisons}

We now present some results of experiments made
with the algorithms described in the previous section, together
with tests performed with a subspace iteration
method (\cite{js}, \cite{sj}) and an Arnoldi method (\cite{tb}, \cite{cb}
\cite{le}, \cite{as}).
The test pencil $(J,L)$ is taken from a transient
stability model of the South-Southeast interconnected Brazilian power system:
$J$ is the Jacobian matrix at an operating point and is of order 3156
and $L$ is a diagonal matrix with elements 1 or 0, corresponding respectively
to state and algebraic variables, with rank equal to 790.
Only $0.14\%$ of the elements of $J$ are nonzero and
its sparseness pattern is given in Figure 1.
This pencil has exactly four eigenvalues
with positive real part: $0.1814 \pm i\, 4.8323$, $0.0233$ and $0.0004$.
The first two correspond
to genuine unstable modes of the system. The remaining two
are related to two redundant states of the system \cite{ku}:
they would have been zero if there were no roundoff errors
in the generation of the Jacobian matrix.

\begin{figure}[t]
\centerline{\psfig{file=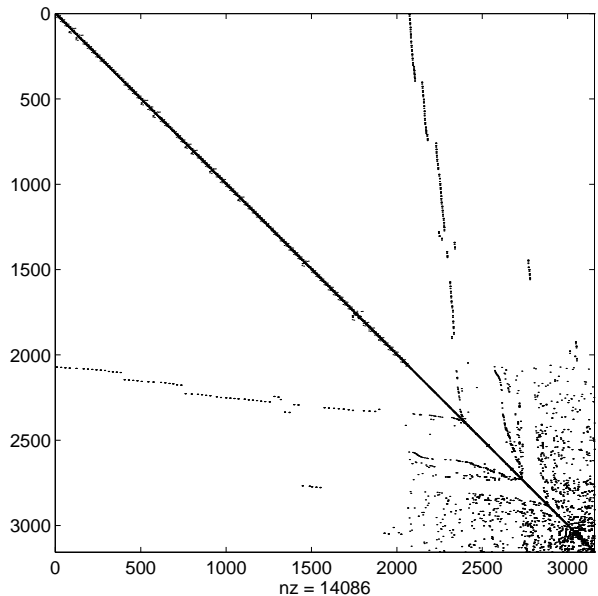,height=8cm}}
\caption{Sparse Pattern of the Test Matrix}
\end{figure}

The same pencil has been used in \cite{cn}, where the authors
report results of a parallelization
of the lopsided simultaneous iteration algorithm \cite{sj}.
Their strategy to find unstable modes is to perform
shift-invert subspace iterations
with initial shifts given on the imaginary axis.
We implemented here a sequential version of this algorithm.
As initial vectors we took $Z^{(0)}=LZ^{(0)}=\pmatrix{X^{(0)}\cr 0}$,
where $X^{(0)}$ is a matrix with (up to) eight first column vectors
of the $790\times 790$ Fourier matrix
$F(j,k)=e^{j.k.\frac{2\pi i}{790}}$.
After each four shift-invert iterations followed by
a normalization of the vectors a Rayleigh-Ritz acceleration was done.
That is,
we first calculated a spectral decomposition of $B=G^{-1}H$
(\cite{la}, \cite{sm}),
where $G=(LZ^{(4k)})^H (LZ^{(4k)})$ and $H=(LZ^{(4k)})^H W$, with
$W=L\, (J-aL)^{-1}L\, Z^{(4k)}$,
such that the Ritz values were ordered according to decreasing absolute values.
Then $W$ was multiplied by the matrix of Ritz vectors and
the resulting column vectors, after
normalization, undertook another cycle of four iterations.
Convergence was achieved when the $\infty - norm$ of the difference
between corresponding vectors of $LZ^{(4k)}$ and $LZ^{(4k+1)}$
was less than a tolerance, taken to be  $10^{-5}$.
A deflation technique was also implemented for these tests:
aside from the choice of initial vectors and the deflation technique, this
is the algorithm used in \cite{cn}.
Some results obtained with these iterations
on SUN SPARC workstations are displayed in Table 1,
where {\it iter} and {\it prod} mean respectively
the number of Rayleigh-Ritz accelerations and
the number of matrix-vector multiplications being performed.
Notice that matrix-vector multiplication consists of
backward and forward substitutions
after computing  an ${\cal LU}$ decomposition of $(J-a L)$,
which is done only once for each $a$.

\begin{table}

\begin{center}
\begin{tabular}{|r|r|r|}  \hline
converged value & iter & prod  \\
\hline
-0.1164+3.2018i
&   3  & 96 \\
-0.0925+3.9827i
&   6  & 180 \\
-0.9573+2.1594i
&   8  & 228 \\
-0.4803+1.7632i
&   8  & 228 \\
\hline
\end{tabular} \hskip 1pc
\begin{tabular}{|r|r|r|}  \hline
converged value & iter & prod  \\
\hline
-0.0925+3.9827i
&   2  & 64 \\
-0.1164+3.2018i
&   7  & 204 \\
0.1814+4.8323i
&   7  & 204 \\
-0.5911+4.6935i
&   8  & 224 \\
\hline
\end{tabular}
\end{center}

\caption{Simultaneous shift-invert method: $a=0+3i$ and $a=0+4i$}

\end{table}

ARNCHEB performs an incomplete
Arnoldi method combined with an acceleration technique
using Chebyshev polynomials
\cite{tb}, \cite{cb}.
We have used its Reverse Communication interface
to calculate matrix-vector multiplications by 
$C_{\sigma}=I + 2\, Re\,\sigma \,L\,(J-\sigma L)^{-1}$.
As seen before, since $L=\pmatrix{I&0\cr 0&0}$,
we avoid the spurious eigenvalue 1
with this Cayley transform
by considering only the first 790 coordinates of the vectors, that is,
by solving systems $(J-\sigma L)w=Lz$ and taking only
the first 790 coordinates of $w$.
In the tests,
the dimension of the Krylov space was taken to be 54
and the number of requested eigenvalues,  4.
Some results are in Table 2, where
${\cal O(r)}$ is the order of the residual
$||(C_{\sigma}-\lambda I) v||_2$, with $||v||_2 = 1$,
and {\it iter} is the number of Arnoldi steps.

\begin{table}

\begin{center}
\begin{tabular}{|r|r|}  \hline
converged value & ${\cal O(r)}$ \\
\hline
0.1814+4.8323i & $10^{-10}$\\
0.1814-4.8323i & $10^{-10}$\\
0.0233+0.0000i & $10^{-08}$\\
0.0004+0.0000i & $10^{-07}$\\
\hline
\hline
\multicolumn{1}{|c|}{iter} & \multicolumn{1}{c|}{prod} \\
\hline
\multicolumn{1}{|c|}{12}  &  \multicolumn{1}{c|}{856} \\
\hline
\end{tabular} \hskip 2pc
\begin{tabular}{|r|r|}  \hline
converged value & ${\cal O(r)}$ \\
\hline
0.1814+4.8323i & $10^{-11}$\\
0.1814-4.8323i & $10^{-11}$\\
0.0233+0.0000i & $10^{-08}$\\
0.0004+0.0000i & $10^{-08}$\\
\hline
\hline
\multicolumn{1}{|c|}{iter} & \multicolumn{1}{c|}{prod} \\
\hline
\multicolumn{1}{|c|}{31}  &  \multicolumn{1}{c|}{2107} \\
\hline
\end{tabular}
\end{center}

\caption{ARNCHEB: $\sigma = 4.0$ and $\sigma = 6.0$}

\end{table}

The same sort of experiment was carried out with ARPACK.
The program dndrv1.f was rewritten to include the same routines which
solved the systems in the tests with ARNCHEB.
Here the dimension of the Krylov space was taken to be 20 and
also 4 eigenvalues were requested to converge.
Some of the results are shown in Table 3.

\begin{table}

\begin{center}
\begin{tabular}{|r|r|}  \hline
converged value & ${\cal O(r)}$ \\
\hline
0.1814+4.8323i & $10^{-13}$\\
0.1814-4.8323i & $10^{-13}$\\
0.0233+0.0000i & $10^{-14}$\\
0.0004+0.0000i & $10^{-13}$\\
\hline
\hline
\multicolumn{1}{|c|}{iter} & \multicolumn{1}{c|}{prod} \\
\hline
\multicolumn{1}{|c|}{301}  &  \multicolumn{1}{c|}{3349} \\
\hline
\end{tabular} \hskip 2pc
\begin{tabular}{|r|r|}  \hline
converged value & ${\cal O(r)}$ \\
\hline
0.1814+4.8323i & $10^{-14}$\\
0.1814-4.8323i & $10^{-14}$\\
0.0233+0.0000i & $10^{-14}$\\
0.0004+0.0000i & $10^{-13}$\\
\hline
\hline
\multicolumn{1}{|c|}{iter} & \multicolumn{1}{c|}{prod} \\
\hline
\multicolumn{1}{|c|}{301}  &  \multicolumn{1}{c|}{3679} \\
\hline
\end{tabular}
\end{center}

\caption{ARPACK: $\sigma = 4.0$ and $\sigma = 6.0$}

\end{table}

The advantage of these methods is that they require only one
factorization of the pencil in upper and lower triangular factors.
The disadvantage is the presence of parameters which need to be adjusted,
like the dimension of the Krylov space
and the convergence criterion.
ARNCHEB worked well when this dimension was large compared to the number of
desired eigenvalues (ten times, e.g.). The tolerance used was 1.0d-11.
The opposite ocurred with ARPACK: it worked well when
the dimension of the Krylov space was between four and six times the number
of desired eigenvalues. The tolerance employed was 0.d0:
when it was changed to 1.0d-16,
the convergence for the two positive real eigenvalues was not achieved.

\begin{table}

\begin{center}
\begin{tabular}{|r|r|r|}  \hline
converged value & iter & ${\cal O(r)}$ \\
\hline
0.0004+0.0000i &  6 (2)&$10^{-10}$\\
-0.6223+0.9649i&  11 (3)&$10^{-10}$\\
-0.4803+1.7632i&  8 (2)&$10^{-07}$ \\
-0.0925+3.9827i & 6 (2)& $10^{-11}$\\
-0.1351+6.8974i & 8 (2)& $10^{-07}$\\
-1.5684+12.593i & 7 (2)& $10^{-08}$\\
\hline
\end{tabular} \hskip 2pc
\begin{tabular}{|r|r|r|}  \hline
converged value & iter & ${\cal O(r)}$ \\
\hline
0.0004+0.0000i & 6 (2)& $10^{-10}$\\
0.0233+0.0000i & 10 (3)& $10^{-11}$\\
0.1814+4.8323i & 10 (3)& $10^{-08}$\\
-0.0925+3.9827i & 6 (2)& $10^{-09}$\\
-0.1351+6.8974i & 9 (3)& $10^{-11}$\\
-0.1144+10.617i & 8 (2)& $10^{-06}$\\
\hline
\end{tabular}
\end{center}

\caption{Algorithm I - $\sigma = 4.0$: $p=0$ and $p=40$}

\end{table}

The experiments with  algorithms I and II were carried out
by taking only the first four column vectors
of the Fourier matrix of order 790 as initial vectors.
We chose initial shifts on the upper half of the unit circle,
$\mu_k =  \displaystyle{e^{\frac{k\pi i}{6}}},
\; 0\leq k < 6$  ($k=6$ yields ($C_{\sigma}^{-1} - I$), which is
singular).
For each initial shift, the process stops when,
for some $k$, $|x_k^{(j)} - x_k^{(j-1)}|< tol$. The tolerance $tol$ was
taken as $10^{-4}$.
The number $t$ of iterations performed before shift
update was 4.
The following tables contain the results of these tests
performed in SUN SPARC workstations.
The first column indicates the eigenvalue reached by the iteration.
Also, {\it iter} is the number of shift-invert iterations,
the number of ${\cal LU}$ factorizations
appears between parenthesis, and
${\cal O(r)}$ is the order of the residual
$$
||(C_{\sigma}^{-1} - \mu\, I) x||_2, \quad ||x||_2=1.
$$
Preconditioning is performed by multiplying initial vectors $p$ times
by $C_{\sigma}$, followed by normalization:
$p=0$ means no preconditioning. If
$p\ne 0$, multiplication by $C_{\sigma}$ requires an additional
${\cal LU}$ decomposition.
In order to know how many matrix-vector multiplications were performed
in the tests, multiply {\it iter} by 4 and add $4p$.

\begin{table}

\begin{center}
\begin{tabular}{|r|r|r|}  \hline
$\lambda$ & iter & ${\cal O(r)}$ \\
\hline
0.0004+0.0000i & 6 (2)&$10^{-10}$\\
-0.4803+1.7632i& 7 (2)&$10^{-07}$\\
-0.1164+3.2018i& 8 (2)&$10^{-07}$\\
-0.1764+6.1231i& 7 (2)&$10^{-08}$\\
-0.1144+10.617i & 5 (2)& $10^{-10}$\\
-1.5684+12.593i & 11 (3)& $10^{-09}$\\
\hline
\end{tabular} \hskip 2pc
\begin{tabular}{|r|r|r|}  \hline
converged value & iter & ${\cal O(r)}$ \\
\hline
0.0004+0.0000i & 6 (2)& $10^{-10}$\\
0.0233+0.0000i & 10 (3)& $10^{-08}$\\
0.1814+4.8323i & 10 (3)& $10^{-08}$\\
-0.1764+6.1231i&  8 (2)& $10^{-07}$\\
-0.1144+10.617i & 4 (1)& $10^{-06}$\\
-101.95+0.0000i & 18 (5)& $10^{-09}$\\
\hline
\end{tabular}
\end{center}

\caption{Algorithm I - $\sigma = 6.0$: $p=0$ and $p=40$}

\end{table}

For the algorithm II, $\sigma$ was chosen to be 4.8334, the modulus of
the unstable eigenvalues $0.1814\pm 4.8323i$. From the previous section,
this is the real value for $\sigma$ that maximizes the absolute value
of the corresponding eigenvalues of the Cayley transform $C_{\sigma}$.
Thus, preconditioning of the initial vectors with this $\sigma$
should make these eigenvalues easier to detect. Indeed,
in Table 6 we can see that one of them was identified from two
consecutive initial shifts. In the same table are listed the results
when the procedure of inhibiting the last found eigenvalue was applied
(for the first initial shift, there is nothing to inhibit).
In Table 7 a better performance of the inhibiting procedure can be seen:
the convergence to the two unstable eigenvalues was achieved.

The ${\cal LU}$ decomposition of $(J-aL)$ is about 150 times slower than the
resolution of the corresponding systems on a SUN SPARC 4 workstation
(typical runtimes were 7.05078 and 4.29688e-02, respectively).
Thus,  Algorithms I and II had a performance comparable to the others
in respect to time and accuracy. For instance,
from Table 4, we see that by preconditioning the vectors we obtained
three unstable modes and three stable ones (two of these are low damped),
after 196+160 matrix-vector multiplications (indeed, backward and forward
substitutions) and 15+1 factorizations.

\begin{table}

\begin{center}
\begin{tabular}{|r|r|r|}  \hline
$\lambda$ & iter & ${\cal O(r)}$ \\
\hline
0.0004+0.0000i & 6 (2)&$10^{-10}$\\
0.0233+0.0000i & 10 (3)& $10^{-09}$\\
0.1814+4.8323i & 12 (3)& $10^{-06}$\\
0.1814+4.8323i &  6 (2)& $10^{-08}$\\
-0.1144+10.617i & 9 (3)& $10^{-12}$\\
-200.38+0.0000i & 15(4)& $10^{-09}$\\
\hline
\end{tabular} \hskip 2pc
\begin{tabular}{|r|r|r|}  \hline
converged value & iter & ${\cal O(r)}$ \\
\hline
0.0004+0.0000i & 6 (2)& $10^{-10}$\\
-0.1764+6.1231i & 12(3) & $10^{-07}$\\
0.0233+0.0000i & 10 (3)& $10^{-09}$\\
0.1814+4.8323i &  6 (2)& $10^{-07}$\\
-0.1144+10.617i & 8 (2)& $10^{-08}$\\
-200.38+0.0000i & 11(3)& $10^{-12}$\\
\hline
\end{tabular}
\end{center}

\caption{Algorithm I and II - $\sigma = 4.8334$, $p=40$}

\end{table}

\begin{table}

\begin{center}
\begin{tabular}{|r|r|r|}  \hline
$\lambda$ & iter & ${\cal O(r)}$ \\
\hline
0.0004+0.0000i & 6 (2)&$10^{-10}$\\
0.0233+0.0000i & 8 (2)& $10^{-06}$\\
0.1814+4.8323i & 7 (2)& $10^{-06}$\\
0.1814+4.8323i &  6 (2)& $10^{-10}$\\
-0.1144+10.617i & 13(4)& $10^{-12}$\\
-0.1144+10.617i &  8(2)& $10^{-07}$\\
\hline
\end{tabular} \hskip 2pc
\begin{tabular}{|r|r|r|}  \hline
converged value & iter & ${\cal O(r)}$ \\
\hline
0.0004+0.0000i & 6 (3)&$10^{-10}$\\
0.1814+4.8323i & 13 (5)& $10^{-09}$\\
0.1814-4.8323i &  9 (4)& $10^{-08}$\\
0.1814+4.8323i & 6 (3)& $10^{-09}$\\
-0.1144+10.617i & 7 (3)& $10^{-07}$\\
0.1814-4.8323i  & 10 (4)& $10^{-07}$\\
\hline
\end{tabular}
\end{center}

\caption{Algorithm I and II - $\sigma = 4.8334$, $p=80$}

\end{table}


\begin{thebibliography}{99}
\bibitem{lh}
L. H. Bezerra, {\it Stability Analysis of Large Scale Power Systems}
(in Portuguese), D. Thesis, Pontif\'{\i}cia Universidade Cat\'olica, Rio
de Janeiro, 1990.
\bibitem{tb}
T. Braconnier, 'The Arnoldi Chebyshev Algorithm for Solving Large Nonsymmetric
Eigenproblems', {\it Technical Rep. TR/PA/93/25}, CERFACS, Toulouse.
\bibitem{cb}
T. Braconnier, V. Fraysse and J.-C. Rioual,
'ARNCHEB Users' Guide: Solution of Large Nonsymmetric or Non Hermitian
Eigenvalue Problems by the Arnoldi-Chebyshev Method',
{\it Technical Rep. TR/PA/97}, CERFACS, Toulouse.
\bibitem{cn}
J. M. Campagnolo, N. Martins, J. L. R. Pereira, L. T. G. Lima, H. J.
C. P. Pinto, and D. M. Falc\~ao, 'Fast Small-Signal Stability
\mbox{Assessment} Using Parallel Processing',
{\it IEEE Trans. on Power Systems}, {\bf PWRS-9(2)}, 949-956, 1994. 
\bibitem{ch}
F. Chatelin and D. Ho, 'Arnoldi-Chebyshev Procedure for Large Scale Nonsymmetric
Matrices', {\it Math. Model. Numer. Anal.}, {\bf 24}, 53-65, 1990.
\bibitem{cg}
K. A. Cliffe, T. J. Garratt and A. Spence, 'Eigenvalues of the Discretized
Navier-Stokes Equation with Application to the Detection of Hopf Bifurcations',
{\it Advances in Computational Mathematics}, {\bf 1(2)}, 337-356, 1993.
\bibitem{ga}
K. A. Cliffe, T. J. Garratt and A. Spence, 'Eigenvalues of Block Matrices
Arising from Problems in Fluid Mechanics', {\it SIAM J. Matrix Anal. Appl.},
{\bf 15(4)}, 1310-1318, 1994. 
\bibitem{ds}
I. S. Duff and J. A. Scott, 'Computing Selected Eigenvalues of Sparse
Unsymmetric Matrices Using Subspace Iteration', {\it ACM Trans.
Math. Softw.}, {\bf 19(2)}, 137-159, 1993. 
\bibitem{go}
G. H. Golub and C. F. Van Loan, {\it Matrix Computations}, 2nd
ed., The Johns Hopkins University Press, Baltimore, 1989.
\bibitem{je}
A. Jennings and J. J. McKeowen, {\it Matrix Computations}, 2nd ed.,
John Wiley and Sons, New York, 1992.
\bibitem{js}
A. Jennings and W. J. Stewart, 'Simultaneous Iteration for Partial
Eigensolution of Real Matrices', {\it J. Inst. Maths. Applics.},
{\bf 15}, 351-361, 1975.
\bibitem{ku} P. Kundur, {\it Power System Stability and Control},
McGraw Hill, New York, 1994.
\bibitem{la}
E. Anderson, Z. Bai, C. Bischof, J. Demmel, J. Dongarra, J. DuCroz,
A. Greenbaum, S. Hammarling, A. McKenney, S. Ostrouchov, and D. Sorensen,
{\it LAPACK Users' Guide, Release 2.0, 2nd ed.},
SIAM Publications, Philadelphia, 1995.
\bibitem{le}
R. B. Lehoucq, D. C. Sorensen and C. Yang, 'ARPACK Users' Guide:
Solution of Large Scale Eigenvalue Problems by Implicitly Restarted
Arnoldi Methods', July 1996.
\bibitem{li}
L. T. G. Lima, L. H. Bezerra, C. Tomei and N. Martins, 'New
Methods for Fast Small Signal Stability Assessment of Large Scale Power
Systems', {\it IEEE Transactions on Power Systems}, {\bf 10(4)},
1979-1985, November 1995.
\bibitem{nm}
N. Martins, 'Efficient Eigenvalue and Frequency Response
Methods Applied to Power System Small-Signal Stability Studies',
{\it IEEE Trans. on Power Systems}, {\bf PWRS-1(1)}, 217-226,
February 1986.
\bibitem{me}
K. Meerbergen and D. Roose, 'Matrix Transformations for Computing Rightmost
Eigenvalues of Large Sparse Non-Symmetric Eigenvalue Problems',
{\it IMA J. Numer. Anal.}, {\bf 16}, 297-346, 1996.
\bibitem{sa}
Y. Saad, 'Chebyshev Acceleration Techniques for Solving Nonsymmetric Eigenvalue
Problems', {\it Math. Comp.}, {\bf 42}, 567-588, 1984.
\bibitem{as}
J. A. Scott, 'An Arnoldi Code for Computing Selected Eigenvalues of Sparse
Real Unsymmetric Matrices', {\it ACM Trans. Math. Softw.},
{\bf 21(4)}, 432-475, 1995.
\bibitem{sm}
B. T. Smith et al., {\it Matrix Eigensystem Routines: EISPACK
Guide}, 2nd ed., Springer Verlag, New York, 1976.
\bibitem{so}
D. C. Sorensen, 'Implicit Application of Polynomial Filters in a K-Step Arnoldi
Method', {\it SIAM J. Matrix Anal. Appl.}, {\bf 13(1)}, 357-385, 1992.
\bibitem{sj}
W. J. Stewart and A. Jennings, 'A Simultaneous Iteration Algorithm for Real
Matrices', {\it ACM Trans. Math. Softw.}, {\bf 7(2)}, 184-198, 1981.
\bibitem{jp}
N. Uchida and T. Nagao,  'A New Eigen-Analysis Method of
Steady-State Stability Studies for Large Power Systems: S-Matrix
Method', {\it IEEE Trans. on Power Systems}, {\bf PWRS-3(2)}, 706-714,
May 1988.
\bibitem{ws}
L. Wang and A. Semlyen, 'Application of Sparse Eigenvalue
Techniques to the Small-Signal Stability Analysis of Large Power
Systems', {\it IEEE Trans. on Power Systems}, {\bf PWRS-5(2)}, 635-642,
May 1990.
\bibitem{wa}
D. S. Watkins, {\it Fundamentals of Matrix Computations},
John Wiley and Sons, New York, 1991.
\end{thebibliography}
\end{document}